\documentclass[DIV=15, bibliography=totoc]{scrartcl}  

\usepackage[a4paper,margin=1.5cm,bottom=2.5cm]{geometry}
\usepackage{amsmath,amssymb}
\usepackage{algorithmicx}
\usepackage{algpseudocode} 
\usepackage{algorithm}

\usepackage{multicol}
\usepackage{nomencl}
\makenomenclature

\renewcommand{\nompreamble}{\begin{multicols}{2}}
\renewcommand{\nompostamble}{\end{multicols}}
\usepackage{float}
\usepackage{xcolor} 
\definecolor{lightblue}{rgb}{0,0.5,1.0}
\definecolor{linkblue}{rgb}{0,0.1,0.6}
\definecolor{citegreen}{rgb}{0,0.4,0.0}%{0.1,0.5,0.4}%{0.125,0.6,0.5}
\definecolor{linkred}{rgb}{0.8,0,0.005}%{0.6,0,0.1}
\definecolor{mailviolet}{rgb}{0.3,0,0.35}%{0.6,0,0.1}
\definecolor{tumblue}{rgb}{0,0.396,0.741}
\definecolor{darkgreen}{rgb}{0,0.4,0} 
\definecolor{darkbrown}{rgb}{0.5, 0.396, 0.09}

\usepackage[hyperindex,colorlinks=true,linkcolor=linkblue,citecolor=citegreen,urlcolor=mailviolet,filecolor=linkred]{hyperref}

\usepackage{caption, subcaption}
\usepackage{fancyhdr}

\usepackage{siunitx}
\usepackage[square,comma,nonamebreak,compress,numbers]{natbib}

\usepackage{soulpos}
\usepackage{ulem}
\usepackage{authblk} % used for affiliations 
\usepackage{float}

\usepackage[inkscapearea=page]{svg}
\usepackage{amsmath}
\usepackage{amssymb}
\usepackage[noabbrev]{cleveref}
\usepackage{upgreek}
\usepackage{adjustbox}
\usepackage{color}
\usepackage{enumitem}
\usepackage{todonotes}

\usepackage{pgfplots}
\usepackage{pgfplotstable}
\usepackage{filecontents}
\usepackage{tikz}
\usepackage{tikzscale}
\usepackage{tikz-3dplot}

\usetikzlibrary{spy}
\usetikzlibrary{intersections}
\usetikzlibrary{arrows,shapes}
\usetikzlibrary{spy}
\usetikzlibrary{backgrounds}  
\usetikzlibrary{decorations}
\usetikzlibrary{decorations.markings}
\usetikzlibrary{positioning}
\usetikzlibrary{patterns}
\usetikzlibrary{calc} 
\usetikzlibrary{quotes} 
\usetikzlibrary{external} % Is activated with the command '\tikzexternalize'.
\usetikzlibrary{matrix}

\pgfplotscreateplotcyclelist{customCycleList}
{%
  {blue, mark=o},
  {red, mark=square},
  {darkgreen, mark=diamond},
  {cyan, mark=diamond},
  {darkbrown, mark=triangle*},
  {black, mark=pentagon},
}

\pgfplotsset{every axis/.append style= {
    cycle list name=customCycleList,
}}

\usepackage{abstract}
\usepackage{placeins}

\usepackage{amssymb, amsmath,bm, amsthm, mathrsfs}

\title{On the radial discretization of finite element spaces in the scaled boundary finite element method}

\author{Alireza Danshyar\thanks{\href{mailto:alireza.daneshyar@tum.de}{\texttt{alireza.daneshyar@tum.de}}, Corresponding author} }
\author{Stefan Kollmannsberger}

\affil[]{Chair of Computational Modeling and Simulation, School of Engineering and Design, Technical University of Munich, Arcisstraße 21, 80333 Munich, Germany}%,Arcisstr. 21, 80333 M\"unchen, Germany}

\newcommand{\publicationDate}{\today}

\date{}
\fancypagestyle{plain}{%
\fancyhf{}%
\vspace{-0.0cm}
\fancyfoot[L]{
\vspace{-0.0cm} 
\footnotesize{\textit{Preprint submitted to Elsevier}\hfill
\publicationDate
\\
}
}
}

\crefname{paragraph}{paragraph}{paragraphs}
\Crefname{paragraph}{Paragraph}{Paragraphs}

\begin{document}
\vspace{-1.5cm} 
\normalem \maketitle  
\normalfont\fontsize{11}{13}\selectfont
\vspace{-1.5cm} \hrule 
\section*{Abstract}
    The scaled boundary finite element method is known for its capability in reproducing highly-detailed solution fields over large subdomains, and also for its flexibility in terms of spatial discretization enabling to define subdomains with arbitrary numbers of boundary elements and hanging nodes. The former, however, is only attainable in those cases where analytical solutions to their differential equations exist. Many others, especially those with oscillatory nature, invoke the use of numerical methods that only provide the response of boundaries. Hence, no information on the inner-subdomain solution fields can be recovered, leaving one of the main asset of the method unused. As a remedy, we propose a new solution scheme by which the interior fields of subdomains can be recovered so that the full potential of the scaled boundary ﬁnite element method can be harvested. No auxiliary variables are introduced to the global algebraic system and the dimensions of the matrices remain intact. Most of the computations are applied prior to the assembly of the global algebraic system, and they are carried out in a completely decoupled manner. This renders the introduced method to be embarrassingly parallelizable. The interior information is not required for the solution process and it can be computed locally for the selected regions and subdomains, also using parallel computing techniques.

    \vspace{0.25cm}
    \noindent\textit{Keywords:} 
    scaled boundary ﬁnite element method; elastodynamics; wave propagation; frequency domain
    \vspace{-0.4cm}

\section{Introduction} \label{sec:introduction}
The scaled boundary finite element method is one of the more recent powerful numerical tools for dealing with problems in the fields of applied science and engineering. It benefits from a unique transformation of geometry that enables defining finite element spaces with arbitrary numbers of boundary elements. According to this transformation, the region occupied by an element, also known as a subdomain, is the set of all points defined by scaling the boundary elements of the subdomain with respect to a scaling center. Thus, it is only necessary to interpolate the field variables along those boundary elements and the field variables along the radial rays that emanate from the scaling center can be kept analytical. By applying the scaled boundary finite element transformation of geometry, the partial differential equations system of the problem is turned into a weak form along the discretized directions while it remains strong along the radial ones. Owing to this so-called \textit{semi-discretized} representation, the scaled boundary finite element method delivers some unique features such as the capability of satisfying radiation conditions at infinity and reproducing the stress singularity. More importantly, its semi-strong form alleviates from the requirement for fines meshes in confronting sharp changes in the field variables, such as those that occur in elastodynamic problems involving high-frequency waves. Nevertheless, there is a trade-off for the effectiveness of this technique in contrast to those using conventional discretization approaches, such as the standard finite element method. The semi-strong form of the transformed equations invokes a preliminary solution by which the contribution of a subdomain in the global algebraic system is determined. This preliminary solution can only be evaluated analytically for very limited cases such as non-oscillatory steady-state problems in elastostatics. Many others, including oscillatory steady-state cases like elastodynamic problems in the frequency domain, and all transient ones must be treated numerically.

The scaled boundary finite element method was first introduced for linear elasticity by Song and Wolf~\cite{song1997scaled} and then extended to cope with various problems, such as heterogeneous solids \cite{reichel2023non}, elastodynamics~\cite{gravenkamp2018scaled, gravenkamp2020mass, gravenkamp2020three, zhang2021massively, qu2022time, zhang2022asynchronous}, nonlinear mechanics~\cite{eisentrager2020sbfem, chasapi2020geometrically, chasapi2020patch, chasapi2021isogeometric}, plates and shells~\cite{wallner2020scaled, li2020efficient, klassen2021isogeometric}, fracture~\cite{song2002semi, yang2006fully, bulling2019high, pramod2019adaptive, zhang2020adaptive, ooi2020polygon}, contact~\cite{ya2021open}, thermoelastic~\cite{iqbal2021development}, ultrasonics~\cite{gravenkamp2017efficient}, acoustics~\cite{liu2019automatic}, to name a few. Approaching the problems of linear elasticity by the method, the system of partial differential equations that governs the displacement field of subdomains is transformed into a boundary value problem that comprises a second-order homogeneous Cauchy--Euler system subjected to Robin boundary conditions~\cite{daneshyar2021general}. This system can be diagonalized by applying a proper linear transformation, leading to a set of uncoupled Cauchy--Euler system. Hence, the analytical solution can be readily obtained and utilized to compute the stiffness matrix of the subdomain. The frequency-domain representation of the elastodynamic problems, on the other hand, forms a set of Bessel differential equations for which an analytical solution cannot be found in general. Transient cases are also transformed into initial-boundary value problems for which numerical solutions are usually employed.

Several strategies have been proposed to substitute analytical solution methods. The earliest and most simple technique was the low-frequency expansion method of Song and Wolf~\cite{song1997scaled} within which the higher-order terms of the dynamic stiffness matrix are dropped, leading to an approximation similar to that of the standard finite element method. As a result, the method looses its superiority over the standard discretization techniques in terms of accuracy. Song and Wolf~\cite{song1998scaled} developed an approximation technique based on an infinite series of matrix manipulation that only yielded acceptable results for low frequencies. Yan et al.~\cite{yann2004coupling} utilized the unit-impulse response functions of Paronesso and Wolf~\cite{paronesso1998recursive} to incorporate unbounded subdomains in three-dimensional dynamic problems. Despite their efforts towards reducing computations associated with matrix manipulation tasks, their solution strategy, however, was not efficient in terms of computational cost. As a result, the improvement of their method in terms of numerical efficiency has been targeted by many authors, among which the works of Genes~\cite{genes2012dynamic} and Schauer et al.~\cite{schauer2012parallel} made compelling progress. In contrast to the infinite series of Song and Wolf~\cite{song1998scaled}, the Frobenius solution of Yang et al.~\cite{yang2007frobenius} could provide acceptable results almost regardless of the frequency of the waves involved, yet it was only applicable for bounded media. Later on, a promising solution strategy, named \textit{continued-fraction method}, was initiated in the work of Song and Bazyar~\cite{song2007boundary}, and matured in a series of studies conducted by Bazyar and Song~\cite{bazyar2008continued}, Song~\cite{song2009scaled}, Birk et al.~\cite{birk2012improved}, and Chen et al.~\cite{chen2014high}. Contrary to the low-frequency expansion of Song and Wolf~\cite{song1997scaled}, the higher-order terms are involved to a desired order in the continued fraction--based methods so that the approximations are improved to the extent that the prevailing frequency content of the problem can be preserved. More recently, a numerical solution strategy base on the shooting approach was presented in a series of studies by Daneshyar and co-authors~\cite{daneshyar2021general, daneshyar2021shooting, daneshyar2023scaled}. Their method treated the two-point boundary-value problem that governs the field variable of the subdomain as a system of coupled initial value problems, integrating the subdomain matrices on the fly.

Despite the numerous research articles on solution strategies for the semi-strong form of the scaled boundary finite element equations, due to the lack of a robust approach, the subject still remains to be explored. The low-frequency expansion method offers no advantage over the standard finite element method in terms of convergence rate as the higher-order terms are dropped from the dynamic stiffness matrix of subdomains. Arranging the scaled boundary finite element method in terms of the unit-impulse response matrix leads to prohibitive computational cost especially if a time-domain solution is of interest within which convolution integrals are involved~\cite{chen2014high}. The Frobenius solution is only given in the frequency domain, it is only derived for bounded subdomains, and its accuracy drops as the frequency increases. The well-established continued-fraction method is formulated for a transformed nonlinear first-order ordinary differential equations system that is only valid on the boundaries of subdomains. Hence, it provides no information on the field variables inside. The same limitation holds true for the shooting approach. Although it uses the original second-order differential equations system of the subdomains to integrate the dynamic stiffness matrix, it performs the integration on the fly, making it impossible to recover the interior field variables. This is a significant drawback since the scaled boundary finite element method is known for its capability in reproducing highly-detailed solution fields within large subdomains. If those data cannot be recovered, the advantage of the method over the conventional ones for the problems whose distribution of field variables is of interest would be questionable. Hence, we propose an alternative solution method here to remove this limitation. To this end, we discretize the radial directions of the subdomains so that the original set of second-order differential equations can be used. As a result, in addition to the dynamic stiffness matrix of the subdomain, their interior information can be recovered within a post-processing step with minimal bookkeeping. Following this introduction, the body of this text proceeds as follows. In Section \ref{sec:basics}, a brief overview of the scaled boundary finite element formulation is given. Section \ref{sec:solution} introduces the proposed solution strategy. The verification and sensitivity analysis of the method are presented in Section \ref{sec:results}. Lastly, Section \ref{sec:conclusion} concludes the paper.

\section{Basics} \label{sec:basics}
\begin{figure}
    \centering
    \includegraphics[scale=1.0]{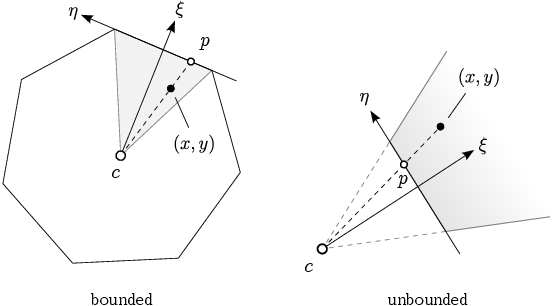}
    \caption{Local coordinates system of the scaled boundary finite element method.}
    \label{fig:local}
\end{figure}

Let us assume an elastic continuous medium subjected to propagating mechanical waves. Its governing equations of motion in the absence of body forces are given in the Voigt notation by
\begin{equation} \label{eq:governing}
    \mathbf L^\intercal\bm\sigma = \rho\ddot{\bm{u}},
\end{equation}
where $\bm\sigma$ is the stress tensor, $\rho$ is the density, $\ddot{\bm{u}}$ is the second time derivative of the displacement, and $\mathbf L$ is the differentiation operator comprising partial derivatives. The scaled boundary finite element method uses a unique coordinate system by which the governing equations of the problem are weakened along all coordinates except one direction, known as the \textit{radial coordinate}. Therefore, whether a two- or three-dimensional elastodynamic problem is of interest, upon applying the scaled boundary transformation of geometry, we obtain a system of second-order ordinary equations. It is worth mentioning that, although the coefficient matrices in a three-dimensional case differ from those of a two-dimensional one, this difference indeed does not affect the solution procedure. Hence, we only address the elastodynamic problems in two-dimensional settings. Accordingly, the operator $\mathbf L$ in two dimensions is given by
\begin{equation}
    \mathbf L =
    \begin{bmatrix}
        \partial_x & 0 \\
        0 & \partial_y \\
        \partial_y & \partial_x
    \end{bmatrix}.
\end{equation}
Employing Hooke's law
\begin{equation}
    \bm\sigma = \mathbf D\bm\varepsilon,
\end{equation}
and the infinitesimal strain theory in which
\begin{equation}
    \bm\varepsilon = \mathbf L\bm u,
\end{equation}
the equation of motion in~\eqref{eq:governing} reads
\begin{equation} \label{eq:matricial}
    \mathbf L^\intercal \mathbf D \mathbf L \bm u = \rho\ddot{\bm{u}}.
\end{equation}
Based on the scaled boundary transformation of geometry, the displacement vector $\bm u$ is defined as
\begin{equation}
    \bm u(\xi,\eta) = \mathbf N(\eta)\bar{\bm u}(\xi) ,
\end{equation}
where $\xi$ and $\eta$ are the local coordinates in the two-dimensional scaled boundary settings. Figure~\ref{fig:local} illustrates these local coordinates in bounded and unbounded subdomains. It should be mentioned that, upon keeping the radial coordinate $\xi$ undiscretized, the problem remains analytical along this coordinate, enabling the representation of unbounded subdomains.

\begin{figure}
    \centering
    \includegraphics[scale=1.0]{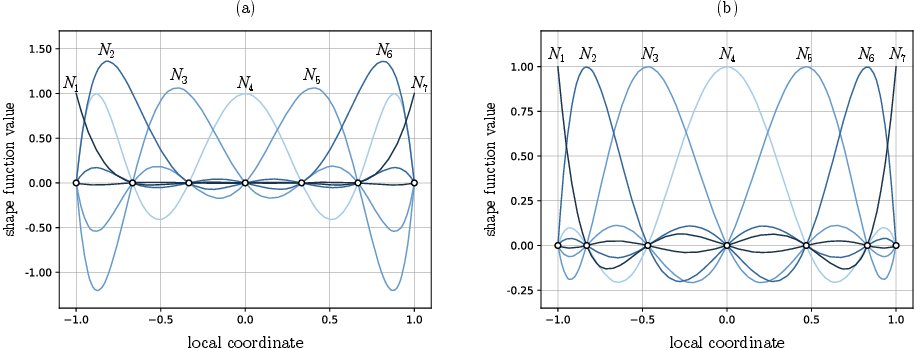}
    \caption{Shape functions of degree six for: (a) evenly-space points, and (b) Gauss--Lobatto quadrature points.}
    \label{fig:polynomials}
\end{figure}

Depending on the choice of interpolation functions, the arrays of $\mathbf N(\eta)$, which is a matrix containing one-dimensional shape functions along the circumferential coordinate $\eta$, are defined. Since, from the theoretical perspective, the field variables are exact along the radial coordinate $\xi$, it is sound to use high-order interpolation functions along the discretized coordinate $\eta$. Adhering to this strategy, the strength of the method to reproduce accurate results in the scaled boundary finite element sub-spaces can be deployed. Regarding the use of higher-order shape functions in the scaled boundary finite element method, Vu and Deeks~\cite{vu2006use} asserted that the so-called $p$-refinement, which corresponds to the degree elevation of interpolation functions, is more effective on the convergence rate rather than the refinement through decreasing the size of sub-spaces, also known as $h$-refinement. They, however, emphasized that, in the case of using evenly-spaced nodes, Lagrange polynomials of relatively high degrees overshoot near the extremities. Furthermore, this elevates the condition number of coefficient matrices and leads to poor approximations. Alternatively, the nodes can be defined at the Gauss--Lobatto quadrature points, leading to well-conditioned coefficient matrices regardless of the degree of polynomials. These points can be given as the roots
\begin{equation}
    (1-\eta^2)P^\prime_n(\eta)=0,
\end{equation}
where
\begin{equation}
    P_n(\eta)=\frac{1}{2^nn!}\frac{\mathrm{d}^n}{\mathrm{d}\eta^n}(\eta^2-1)^n
\end{equation}
is the Rodrigues representation of the Legendre polynomial of degree $n$, and $\eta$ is the local coordinate ranging over $[-1,+1]$, which is the interval bounded between the extremities~\cite{riley1999mathematical}. Defining the shape functions $N_i(\eta)$ using the resulting functions, which are known as the Gauss--Lobatto--Legendre polynomials, we arrive at
\begin{equation}
    N_i(\eta) = \prod_{\substack{j=1 \\ j\ne i}}^{n+1}{\frac{\eta - \eta_j}{\eta_i - \eta_j}},
\end{equation}
wherein $\eta_i$ denotes the $i$-th Gauss--Lobatto quadrature point. To compare the choices of the conventional evenly-spaced points and the Gauss--Lobatto quadrature points, the shape functions of degree six for both cases are plotted next to each other in Figure~\ref{fig:polynomials}. As it is evident, the Gauss--Lobatto--Legendre polynomials do not overshoot near extremities and range over a limited interval.

With the shape functions at hand, the projection of an arbitrary point on its respective boundary of the subdomain reads
\begin{subequations}
    \begin{gather}
        x_p(\eta) = \sum_i{N_i(\eta)} x_i  \\
        y_p(\eta) = \sum_i{N_i(\eta)} y_i .
    \end{gather}
\end{subequations}
where $(x_p,y_p)$ and $(x_i,y_i)$ are the coordinates of the projection point and boundary nodes of the simplex, respectively~\cite{wolf2003scaled, song2004matrix}. Now, that arbitrary point is given by means of the scaled boundary transformation of geometry as~\cite{song2018scaled}
\begin{subequations} \label{eq:transformation}
    \begin{gather}
    x(\xi,\eta) = x_c + \xi(x_p(\eta)-x_c) \\
    y(\xi,\eta) = y_c + \xi(y_p(\eta)-y_c)
    \end{gather}
\end{subequations}
wherein $(x_c,y_c)$ is the coordinates of the scaling center, assigned uniquely to each subdomain (see Figure~\ref{fig:local}). It is the point that the radial lines of the local coordinate $\xi$ emanate from. This point must be chosen such that any ray originating from it only crosses the boundary of the subdomain once. This requirement implies that star-convex sets are only valid for discretizing the domain of a problem~\cite{ooi2016construction}. It need be emphasized that the choice of the scaling center can affect the solutions. The reason lies in the fact that the resulting system of equations for problems with oscillatory nature is actually a Bessel matrix differential equation system~\cite{song2006dynamic, kausel2019numerical}. Hence, the method can provide exact solutions for circular and spherical cells having the scaling center at their origin. Therefore, the choice of this point has to be made deliberately to achieve the best possible central symmetry.

Due to employing the local coordinates system, the partial derivatives with respect to the global coordinates must be given in terms of the local ones using the chain rule
\begin{subequations}
    \begin{gather}
        \frac{\partial}{\partial \xi}=\frac{\partial}{\partial x}\frac{\partial x}{\partial \xi}+\frac{\partial}{\partial y}\frac{\partial y}{\partial \xi}, \\
        \frac{\partial}{\partial \eta}=\frac{\partial}{\partial x}\frac{\partial x}{\partial \eta}+\frac{\partial}{\partial y}\frac{\partial y}{\partial \eta},
    \end{gather}
\end{subequations}
or, equivalently,
\begin{equation}
    \begin{bmatrix}
        \partial_{,\xi}  \\
        \partial_{,\eta} 
    \end{bmatrix} =
    \begin{bmatrix}
        x_{,\xi}  & y_{,\xi}  \\
        x_{,\eta} & y_{,\eta} \\
    \end{bmatrix}
    \begin{bmatrix}
        \partial_{,x}  \\
        \partial_{,y} 
    \end{bmatrix} .
\end{equation}
Applying the scaled boundary transformation of geometry~\eqref{eq:transformation}, we arrive at
\begin{equation}
    \begin{bmatrix}
        \partial_{,\xi}  \\
        \partial_{,\eta} 
    \end{bmatrix} =
    \begin{bmatrix}
        1 & 0 \\
        0 & \xi
    \end{bmatrix} 
    \begin{bmatrix}
        x_p-x_c    & y_p-y_c \\
        x_{p,\eta} & y_{p,\eta}
    \end{bmatrix} 
    \begin{bmatrix}
        \partial_{,x}  \\
        \partial_{,y} 
    \end{bmatrix} .
\end{equation}
By rewriting the above chain rule in its inverse form we have
\begin{equation}
    \begin{bmatrix}
        \partial_{,x}  \\
        \partial_{,y} 
    \end{bmatrix} = \frac{1}{\xi}
    \begin{bmatrix}
        j_{11} & j_{12} \\
        j_{21} & j_{22}
    \end{bmatrix}
    \begin{bmatrix}
        \xi\partial_{,\xi}  \\
        \partial_{,\eta} 
    \end{bmatrix} ,
\end{equation}
wherein
\begin{subequations}
    \begin{gather}
        j_{11} =  \frac{1}{|\mathbf J|}y_{p,\eta} , \\
        j_{12} = -\frac{1}{|\mathbf J|}(y_p-y_c) , \\
        j_{21} = -\frac{1}{|\mathbf J|}x_{p,\eta} , \\
        j_{22} =  \frac{1}{|\mathbf J|}(x_p-x_c),
    \end{gather}
\end{subequations}
and
\begin{equation}
    \mathbf J =
    \begin{bmatrix}
        x_p-x_c    & y_p-y_c \\
        x_{p,\eta} & y_{p,\eta}
    \end{bmatrix}
\end{equation}
is the Jacobian matrix on the corresponding boundary. Now, the differentiation operator $\mathbf L$ is defined as
\begin{equation}
    \mathbf L = \mathbf b_1 \partial_{,\xi} + \frac{1}{\xi}\mathbf b_2 \partial_{,\eta} ,
\end{equation}
where
\begin{subequations}
    \begin{gather}
    \mathbf b_1 =
    \begin{bmatrix}
        j_{11} & 0 \\
        0 & j_{21} \\
        j_{21} & j_{11}
    \end{bmatrix} , \\
    \mathbf b_2 =
    \begin{bmatrix}
        j_{12} & 0 \\
        0 & j_{22} \\
        j_{22} & j_{12}
    \end{bmatrix} .
\end{gather}
\end{subequations}
Finally, the strain tensor is given by
\begin{equation}
    \bm\varepsilon = \mathbf B_1\bar{\bm u}_{,\xi} + \frac{1}{\xi}\mathbf B_2\bar{\bm u},
\end{equation}
wherein
\begin{subequations}
    \begin{align}
        \mathbf B_1 & = \mathbf b_1 \mathbf N, \\
        \mathbf B_2 & = \mathbf b_2 \mathbf N_{,\eta}.
    \end{align}
\end{subequations}

Upon applying the Galerkin method of weighted residuals on the system of partial differential equations in~\eqref{eq:matricial}, the scaled boundary finite element equation of elastodynamic problems in two-dimensional settings is obtained as~\cite{song1997scaled}
\begin{equation}
    \xi^2 \mathbf E_0 \bar{\bm u}_{,\xi\xi} + \xi (\mathbf E_0 + \mathbf E_1^\intercal - \mathbf E_1)\bar{\bm u}_{,\xi} - \mathbf E_2\bar{\bm u} - \xi^2\mathbf M_0 \ddot{\bar{\bm u}} = \mathbf 0,
\end{equation}
with the internal force vector
\begin{equation}
    \bar{\bm q} = \xi \mathbf E_0\bar{\bm u}_{,\xi} + \mathbf E_1^\intercal \bar{\bm u},
\end{equation}
where
\begin{subequations}
    \begin{align}
        & \mathbf E_0 = \int_\eta{\mathbf B_1^\intercal \mathbf D \mathbf B_1|\mathbf J|{d}\eta}, \\
        & \mathbf E_1 = \int_\eta{\mathbf B_2^\intercal \mathbf D \mathbf B_1|\mathbf J|{d}\eta}, \\
        & \mathbf E_2 = \int_\eta{\mathbf B_2^\intercal \mathbf D \mathbf B_2|\mathbf J|{d}\eta}, \\
        & \mathbf M_0 = \int_\eta{\mathbf N^\intercal   \rho      \mathbf N  |\mathbf J|{d}\eta},
    \end{align}
\end{subequations}
are some coefficient matrices assembled for each subdomain by iterating over quadrature points on its one-dimensional boundary elements.

To obtain the frequency-domain representation of the resulting initial-boundary value problem, the exponential Fourier transform pair
\begin{subequations}
    \begin{align}
        \hat{\mathscr{F}}(\omega) &= \int_{-\infty}^{+\infty}{\mathscr{F}(t)\exp(-i\omega t)dt}, \\
        \mathscr{F}(t) &= \int_{-\infty}^{+\infty}{\hat{\mathscr{F}}(\omega)\exp(i\omega t)d\omega},
    \end{align}
\end{subequations}
wherein $i$ is the imaginary unit and $\omega$ is the angular frequency, must be applied. Performing the transformation, we arrive at the system of second-order ordinary differential equations
\begin{equation} \label{eq:second-order}
    \xi^2 \mathbf E_0 \hat{\bm u}_{,\xi\xi} + \xi (\mathbf E_0 + \mathbf E_1^\intercal - \mathbf E_1)\hat{\bm u}_{,\xi} - (\mathbf E_2 - \xi^2 \omega^2\mathbf M_0)\hat{\bm u} = \mathbf 0,
\end{equation}
with the internal force vector
\begin{equation}
    \hat{\bm q} = \xi \mathbf E_0\hat{\bm u}_{,\xi} + \mathbf E_1^\intercal \hat{\bm u},
\end{equation}
wherein $\hat{\bm u}$ and $\hat{\bm q}$ are functions of $\xi$ and $\omega$.

The most common strategy to approach the above boundary value problem is to use the continued fraction--based methods~\cite{song2007boundary, bazyar2008continued, song2009scaled, birk2012improved, chen2014high}. Defining the dynamic stiffness matrix $\mathbf S$ on the boundary $\xi=1$ such that
\begin{equation}
    \hat{\bm q} = \mathbf S(\omega) \hat{\bm u},
\end{equation}
those methods transform the system of second-order ordinary differential equations in~\eqref{eq:second-order} into the so-called scaled boundary finite element equation in dynamic stiffness~\cite{song2018scaled}
\begin{equation} \label{eq:first-order}
    (\mathbf S - \mathbf E_1)\mathbf E_0^{-1}(\mathbf S - \mathbf E_1^\intercal ) - \mathbf E_2 + \omega\mathbf S_{,\omega} + \omega^2\mathbf M_0 = \bm 0 ,
\end{equation}
which is a set of nonlinear first-order ordinary differential equations for the dynamic stiffness matrix $\mathbf S$ in terms of the independent variable $\omega$. The continued-fraction solution typically proceeds by setting
\begin{equation}
    \lambda = -\omega^2,
\end{equation}
and assuming that the dynamic stiffness matrix is a summation of constant, linear, and higher-order terms as follows
\begin{equation}
    \mathbf S(\lambda) = \mathbf K + \lambda\mathbf M - \lambda^2[\mathbf S^{(1)}(\lambda)]^{-1}
\end{equation}
wherein $\mathbf S^{(1)}(\lambda)$ comprises the residual terms, which can be obtained by means of the recursive expression
\begin{equation}
    \mathbf S^{(i)}(\lambda) = \mathbf S_0^{(i)} + \lambda\mathbf S_1^{(i)} - \lambda^2[\mathbf S^{(i+1)}(\lambda)]^{-1}.
\end{equation}
Starting from the highest predefined order $n$, which is chosen based on the frequency content of the problem and the longest radial distance from the scaling center, the contribution of the higher-order residual terms is assumed to be negligible, so that
\begin{equation}
    \lambda^2[\mathbf S^{(n+1)}(\lambda)]^{-1} = \bm 0 .
\end{equation}
Now the matrices $\mathbf S_0^{(i)}$ and $\mathbf S_0^{(i)}$ can be computed by plugging the residual terms $\mathbf S^{(i)}(\lambda)$ in the differential equations system~\eqref{eq:first-order}, from the highest to the lowest order. Therefore, the dynamic stiffness matrix $\mathbf S$ is computed through a recursive matrix inversion procedure. Note that the resulting dynamic stiffness matrix is approximate in the sense that the higher the order of residual terms, the more accurate the dynamic stiffness matrix $\mathbf S$. In addition, regardless of how accurate it is, the dynamic stiffness matrix is only valid on the outer boundary of the subdomain. It is because of the assumption made to enable transforming the linear second-order differential equations in~\eqref{eq:second-order} into the nonlinear first-order differential equations in~\eqref{eq:first-order}. As a result, it provides no information on the field variables inside the subdomain, thus rendering the recovery of the internal solution impossible. 

In order to reproduce the filed variables at arbitrary radial distances, instead of solving the scaled boundary finite element equation in dynamic stiffness, the original differential equations system, which is given in terms of the independent variable $\xi$, need be integrated. However, this system is actually a two-point boundary value problem that exploits introducing auxiliary internal unknowns if a numerical discretization method is opted for. As a result, the number of unknowns grows exponentially and the dynamic stiffness matrix of the system becomes a large matrix of matrices. This two-point boundary value problem, however, can be treated as a set of coupled initial value problems. As a result, self-starting numerical integration methods such as those belonging to the Runge--Kutta family can be applied. This technique is the basis of the shooting approach presented by Daneshyar et al.~\cite{daneshyar2021shooting}. However, the dynamic stiffness matrix is integrated on the fly, hence making it impossible to recover the information inside.

\section{Solution} \label{sec:solution}
We demonstrate our solution procedure via the central differencing scheme. However, it can be implemented by means of any numerical integration technique tailored for two-point boundary value problems. Care must be taken regarding the bandwidth of the resulting coefficient matrix since violating its tridiagonality leads to numerical overheads that may not compensate for the improvement made on the convergence rate. It also must be emphasized that the method does not exploit assembling a large dynamic stiffness matrix composed of sub-matrices. Therefore, no auxiliary internal variables are involved and the dimensions of the global algebraic system remain constant. We only introduce such large coefficient matrix in the text to demonstrate the solution procedure.

According to the second-order central difference approximation, we have
\begin{subequations}
    \begin{gather}
    \hat{\bm u}_{,\xi}    = \delta^1\hat{\bm u}_i + O(h^2), \\
    \hat{\bm u}_{,\xi\xi} = \delta^2\hat{\bm u}_i + O(h^2), 
    \end{gather}
\end{subequations}
where
\begin{subequations}
    \begin{gather}
    \delta^1\hat{\bm u}_i = \frac{1}{2h} (-\hat{\bm u}_{i-1} +  \hat{\bm u}_{i+1}), \\
    \delta^2\hat{\bm u}_i = \frac{1}{h^2} (\hat{\bm u}_{i-1} - 2\hat{\bm u}_{i} + \hat{\bm u}_{i+1}),
    \end{gather}
\end{subequations}
where $h$ is the point spacing and $\hat{\bm u}_i$ is the frequency-domain representation of the displacement vector at $\xi_i$, which is the $i$-th discretization point. Applying the central difference approximation, the system in~\eqref{eq:second-order} is discritized as
\begin{equation}
    \xi_i^2 \mathbf E_0\delta^2\hat{\bm u}_i  + \xi_i (\mathbf E_0 + \mathbf E_1^\intercal - \mathbf E_1)\delta^1\hat{\bm u}_i - (\mathbf E_2 - \xi_i^2 \omega^2\mathbf M_0)\hat{\bm u}_i = \mathbf 0,
\end{equation}
for $i\in\lbrace 0,1,2,\cdots,n \rbrace$ so that $\xi_0$ and $\xi_n$ denote the inner and outer boundaries of the subdomain. Thus, we arrive at the following $n+1$ equations
\begin{equation}
    \bm\vartheta_i \hat{\bm u}_{i-1} + \bm\varphi_i \hat{\bm u}_i + \bm\psi_i \hat{\bm u}_{i-1} = \mathbf 0,
\end{equation}
where
\begin{subequations}
    \begin{gather}
    \bm\vartheta_i = \frac{1}{h^2}\xi_i^2\mathbf E_0 - \frac{1}{2h}\xi_i (\mathbf E_0 + \mathbf E_1^\intercal - \mathbf E_1), \\
    \bm\varphi_i   = -\frac{2}{h^2}\xi_i^2 \mathbf E_0 - \mathbf E_2 + \xi_i^2\omega^2 \mathbf M_0, \\
    \bm\psi_i      = \frac{1}{h^2}\xi_i^2\mathbf E_0 + \frac{1}{2h}\xi_i (\mathbf E_0 + \mathbf E_1^\intercal - \mathbf E_1).
    \end{gather}
\end{subequations}
The two ghost vectors $\hat{\bm u}_{-1}$ and $\hat{\bm u}_{n+1}$ are by-products of this spatial discritization. Hence, two more equations are needed, which must be defined by applying the central differencing scheme on the boundary conditions. Discretizing the internal force vector $\hat{\bm q}$ gives
\begin{equation}
    \hat{\bm q}_i = \xi_i \mathbf E_0\delta^1\hat{\bm u}_i + \mathbf E_1^\intercal \hat{\bm u}_i.
\end{equation}
Due to the equilibrium, the internal force vector vanishes at the inner boundary $\xi_0$. Thus, we arrive at
\begin{equation} \label{eq:zero-force}
    \xi_0 \mathbf E_0\delta^1\hat{\bm u}_0 + \mathbf E_1^\intercal \hat{\bm u}_0 = \mathbf 0,
\end{equation}
which can be given by
\begin{equation}
    \bm\vartheta_{-1}  \hat{\bm u}_{-1}  + \bm\varphi_{-1}  \hat{\bm u}_0 + \bm\psi_{-1}  \hat{\bm u}_1     = \mathbf 0,
\end{equation}
wherein
\begin{subequations}
    \begin{gather}
    \bm\vartheta_{-1} = -\frac{1}{2h}\xi_0 \mathbf E_0, \\
    \bm\varphi_{-1}   =  \mathbf E_1^\intercal, \\
    \bm\psi_{-1}      =  \frac{1}{2h}\xi_0 \mathbf E_0.
    \end{gather}
\end{subequations}
Similarly, for the open boundary of the subdomain at $\xi_n$ we have
\begin{equation} \label{eq:outer}
    \xi_n\mathbf E_0\delta^1\hat{\bm u}_n + \mathbf E_1^\intercal \hat{\bm u}_n = \pm\bm f,
\end{equation}
where $\bm f$ is the external force vector acting on the global degrees of freedom. Note that the choice of positive or negative sign depends on whether the subdomain is bounded or unbounded. For the latter case, the negative sign applies since the outward normal on the open boundary points toward the scaling center, which is opposite to the local coordinate $\xi$. The boundary condition in~\eqref{eq:outer} can be rewritten as
\begin{equation}
    \bm\vartheta_{n+1} \hat{\bm u}_{n-1} + \bm\varphi_{n+1} \hat{\bm u}_n + \bm\psi_{n+1} \hat{\bm u}_{n+1} = \pm\bm f,
\end{equation}
where
\begin{subequations}
    \begin{gather}
    \bm\vartheta_{n+1} = -\frac{1}{2h}\xi_n \mathbf E_0, \\
    \bm\varphi_{n+1}   =  \mathbf E_1^\intercal, \\
    \bm\psi_{n+1}      =  \frac{1}{2h}\xi_n \mathbf E_0.
    \end{gather}
\end{subequations}
Finally, the following algebraic system of equations is obtained
\begin{equation} \left[
    \begin{array}{lllllll}
        \bm\vartheta_{-1}&\bm\varphi_{-1} &\bm\psi_{-1}  &                  &                  &                &             \\
        \bm\vartheta_{0} &\bm\varphi_{0}  &\bm\psi_{0}   &                  &                  &                &             \\
                         &\bm\vartheta_{1}&\bm\varphi_{1}&\bm\psi_{1}       &                  &                &             \\
                         &                &\ddots        &\ddots            &\ddots            &                &             \\
                         &                &              &\bm\vartheta_{n-1}&\bm\varphi_{n-1}  &\bm\psi_{n-1}   &             \\
                         &                &              &                  &\bm\vartheta_{n}  &\bm\varphi_{n}  &\bm\psi_{n}  \\
                         &                &              &                  &\bm\vartheta_{n+1}&\bm\varphi_{n+1}&\bm\psi_{n+1}
    \end{array} \right] \left[
    \begin{array}{l}
        \hat{\bm u}_{-1} \\ \hat{\bm u}_0 \\ \hat{\bm u}_1 \\ \vdots \\ \hat{\bm u}_{n-1} \\ \hat{\bm u}_{n} \\ \hat{\bm u}_{n+1}
    \end{array} \right] = \pm\left[
        \begin{array}{l}
            \mathbf 0 \\ \mathbf 0 \\ \mathbf 0 \\ \vdots \\ \mathbf 0 \\ \mathbf 0 \\ \bm f
        \end{array} \right],
\end{equation}
wherein the coefficient is a matrix of matrices, and the unknowns and external forces are vectors of vectors.

Due to the presence of $\bm\psi_{-1}$ and $\bm\vartheta_{n+1}$, the coefficient matrix is not tridiagonal. However, by applying the Gauss--Jordan elimination on the first and last equations sets, the rest can be treated similar to a tridiagonal system. As a result, the system need not be assembled explicitly and the coefficient matrix can be transformed to an upper-triangular matrix on the fly using the Thomas algorithm~\cite{patankar2018numerical}.

Following the Gauss--Jordan elimination method, we define the pivot matrix of the second equations set as
\begin{equation}
    \bm\chi_0 = \bm\vartheta_0\bm\vartheta_{-1}^{-1}.
\end{equation}
Now, by pre-multiplying the first equations set with $\bm\chi_0$ and subtracting the result from the second equations set, the sub-matrices of the second set are redefined as a linear combination of the first and second sets, such that
\begin{subequations}
    \begin{gather}
    \bm\vartheta_0 := \bm\vartheta_0 - \bm\chi_0\bm\vartheta_{-1} = \bm{0}, \\
    \bm\varphi_0   := \bm\varphi_0   - \bm\chi_0\bm\varphi_{-1}, \\
    \bm\psi_0      := \bm\psi_0      - \bm\chi_0\bm\psi_{-1}.
    \end{gather}
\end{subequations}
As a result, the sub-diagonal matrix $\bm\vartheta_{0}$ is eliminated. Similarly, by defining
\begin{equation}
    \bm\chi_{n+1} = \bm\vartheta_{n+1}\bm\vartheta_n^{-1},
\end{equation}
the sub-matrices of the last equations set are corrected as
\begin{subequations}
    \begin{gather}
    \bm\vartheta_{n+1} := \bm\vartheta_{n+1} - \bm\chi_{n+1}\bm\vartheta_n = \bm{0}, \\
    \bm\varphi_{n+1}   := \bm\varphi_{n+1}   - \bm\chi_{n+1}\bm\varphi_n, \\
    \bm\psi_{n+1}      := \bm\psi_{n+1}      - \bm\chi_{n+1}\bm\psi_n,
    \end{gather}
\end{subequations}
which leads to the elimination of $\bm\vartheta_{n+1}$. Next, to eliminate the sub-matrices $\bm\vartheta_{i}$ for $i\in\lbrace 0,1,2,\cdots,n \rbrace$, we define
\begin{equation}
    \bm\chi_i = \bm\vartheta_i\bm\varphi^{-1}_{i-1}
\end{equation}
by which the following forward sweep is performed
\begin{subequations}
    \begin{gather}
    \bm\vartheta_i := \bm\vartheta_i - \bm\chi_i\bm\varphi_{i-1} = \bm{0}, \\
    \bm\varphi_i   := \bm\varphi_i   - \bm\chi_i\bm\psi_{i-1}.
    \end{gather}
\end{subequations}
Now, the last two equations sets, which are decoupled from the system, read
\begin{equation} \label{eq:decoupled} \left[
    \begin{array}{ll}
        \bm\varphi_{n}  &\bm\psi_{n}  \\
        \bm\varphi_{n+1}&\bm\psi_{n+1}
    \end{array} \right] \left[
    \begin{array}{l}
        \hat{\bm u}_{n} \\ \hat{\bm u}_{n+1}
    \end{array} \right] = \pm\left[
    \begin{array}{l}
        \mathbf 0 \\ \bm f
    \end{array} \right].
\end{equation}
The first set gives
\begin{equation}
    \hat{\bm u}_{n+1} = -\bm\psi_{n}^{-1}\bm\varphi_{n}\hat{\bm u}_{n}.
\end{equation}
Therefore, eliminating $\hat{\bm u}_{n+1}$ from the second one leads to
\begin{equation}
    (\bm\varphi_{n+1}-\bm\psi_{n+1}\bm\psi^{-1}_{n}\bm\varphi_{n})\hat{\bm u}_{n} = \pm\bm f.
\end{equation}
Considering that $\xi_n$ corresponds to the outer boundary of the subdomain, the vector $\hat{\bm u}_n$ denotes the nodal displacement vector ${\hat{\bm u}_n}={\hat{\bm u}^g}$, i.e., actually the global degrees of freedom. Thus, we arrive at
\begin{equation} \label{eq:algebraic}
    \mathbf S \hat{\bm u}^g = \bm f,
\end{equation}
wherein
\begin{equation}
    \mathbf S = \pm(\bm\varphi_{n+1}-\bm\psi_{n+1}\bm\psi^{-1}_{n}\bm\varphi_{n}).
\end{equation}
Note that the choice of the positive or negative sign must appear on the left-hand side of~\eqref{eq:algebraic} since the external force vector on the right-hand side is assembled by discretizing the applied load and its sign does not depend on the type of subdomains. As a result, depending on whether the subdomain is bounded or unbounded, a positive or negative sign must be applied on its dynamic stiffness matrix. Now, the global algebraic system, which only involves nodal displacements, can be assembled and solved without incorporating the internal auxiliary unknowns. Hence, the system is compact in the sense that the computations made thus far are performed to obtain the dynamic stiffness matrix of the subdomain, and no large matrices need be assembled.

Now, by solving the global algebraic system, the nodal displacement vector $\hat{\bm u}_n$ can be obtained. At this point, for the purpose of post-processing, the internal fields can be recovered in a similar manner. Therefore, by performing the backward sweep
\begin{equation}
    \hat{\bm u}_i = -\bm\varphi_i^{-1}\bm\psi_i\hat{\bm u}_{i+1},
\end{equation}
starting from the known nodal vector $\hat{\bm u}_n$ on the outer boundary of the subdomain, the internal auxiliary unknowns $\hat{\bm u}_{i}$ can be computed.

\subsection{Domain of integration}
The choice of $\xi_n=1$ is obvious as it denotes the outer boundary of the subdomain. However, the choice of the starting point depends on whether the subdomain is bounded or unbounded. In the former case, the theoretical starting point is $\xi_0=0$. However, the differential equations system in~\eqref{eq:second-order} is singular at the scaling center. This singularity can be easily circumvented by satisfying the boundary condition of the starting point within an asymptotic sense. To this end, instead of setting $\xi_0=0$, a value close to zero can be used. Although smaller values comply better with the mathematics involved, care must be taken since the closer to the singular point, the higher the condition number of the matrices (see~\cite{daneshyar2021general} for a comprehensive sensitivity analysis regarding this choice).

The situation is completely different for the unbounded subdomains. They have a clear-cut at their open boundary, but they extend toward infinity from the other side. Numerical discretization, however, can be applied on bounded intervals. Therefore, the unbounded subdomains must be truncated at a finite distance from their open end. However, this numerical treatment does not resolve the issue by itself since no matter how far the end-point is located, once an outgoing wave travels the whole length, it hits the truncated boundary and is then reflected back toward the physical domain. Hence, the zero-force boundary condition defined in~\eqref{eq:zero-force} only corresponds to the radiation condition if the outgoing waves become completely exhausted. Therefore, energy dissipation through material damping is the key ingredient of unbounded subdomains, otherwise, its truncated boundary resembles a traction-free surface that causes reflection. Fortunately, the frequency-domain representation of the elastodynamic problems can be easily augmented by an imaginary part representing the hysteretic damping and extended for viscoelastic materials. To this end, the constitutive matrix of the viscoelastic media reads
\begin{equation}
    \hat{\mathbf D} = (1+2i\zeta)\mathbf D
\end{equation}
wherein $\zeta$ is the frequency-independent damping ratio. Now, by truncating the unbounded subdomains far enough from the scaling center, the radiation condition can be established. It should be mentioned that longer radial distances invoke more discretization points. Hence, by assigning larger damping ratios to the unbounded subdomains, the required radial distance and the computation cost associated with it can be reduced. However, to prevent wave scattering due to sudden changes in the material properties, the damping ratio of the unbounded subdomain at $\xi=1$ must be identical to that of the physical domain. A remedy is to assume a gradual growth of the damping ratio from the open end toward the truncation point. This numerical treatment also allows a non-dissipative domain to meet a dissipative boundary. The reader is referred to the work of Daneshyar et al.~\cite{daneshyar2021shooting, daneshyar2023scaled} for further details.

\section{Results} \label{sec:results}
\begin{figure}
    \centering
    \includegraphics[scale=1.0]{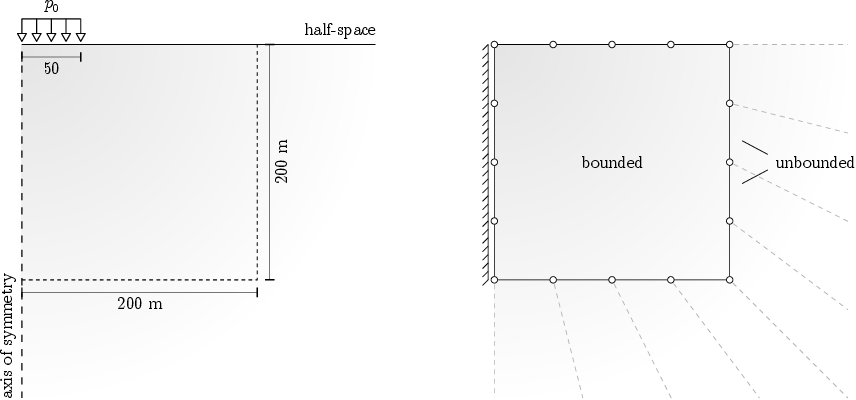}
    \caption{Half-space problem and its scaled boundary finite element mesh.}
    \label{fig:halfspace}
\end{figure}

The analytical solution of a viscoelastic half-space subjected to harmonic strip loading is used to verify the method. The configuration of the problem and its scaled boundary finite element mesh is illustrated in Figure~\ref{fig:halfspace}. Due to symmetry, only the right half of the problem is modeled.

The analytical vertical displacement of the surface in the frequency-domain reads~\cite{miller1954field} 
\begin{equation} \label{eq:Analytical}
    \hat{v}(x)= \frac{2c_s^2p_0}{\pi\mu(1+2i\zeta)^2} \int_0^\infty{
    \frac{\alpha \sin(b\tau)\cos(x\tau)}{\tau(\tau^2+\beta^2)^2 - 4\tau^3\alpha\beta} {d}\tau},
\end{equation}
with
\begin{subequations}
    \begin{align}
        \alpha^2 = \tau^2 - \frac{c_p^2}{(1+2i\zeta)}, \\
        \beta^2  = \tau^2 - \frac{c_s^2}{(1+2i\zeta)},
    \end{align}
\end{subequations}
where $\mu$ is the shear modulus of the half-space, $p_0$ is the magnitude of the strip load, $b$ is the half of the strip load width, and $c_p$ and $c_s$ are respectively the speeds of the pressure and shear waves. The Young’s modulus, Poisson’s ratio, mass density, and damping ratio of the half-space are chosen to be 10 GPa, 0.2, 2500 kg/m$^3$, and 0.05, respectively. The magnitude of the strip load is assumed to be 1.0 GPa, and the problem is solved under plane strain conditions. The integration on the right-hand side should be evaluated by means of a numerical integration method. We used Simpson's second rule, by which
\begin{equation}
    \int_0^\infty{F({x}){d}{x}}\approx\frac{3\Delta{x}}{8}\sum_{i=1}^n
    {\Big( F({x}_{3i-3}) + 3F({x}_{3i-2}) + 3F({x}_{3i-1}) + F({x}_{3i}) \Big)},
\end{equation}
where
\begin{equation}
    \Delta{x} = \frac{{x}_{3n}-x_0}{3n}.
\end{equation}
Note that assigning a sufficiently large value to the upper limit guarantees reaching the asymptotic accuracy. Since the problem is analyzed in the frequency domain, the solutions are combinations of real and imaginary parts. Hence, the time variation of the displacement field is resembled using
\begin{equation}
    \bm u(\phi) = \mathrm{Re}(\hat{\bm u})\cos(\phi) - \mathrm{Im}(\hat{\bm u})\sin(\phi),
\end{equation}
wherein $\phi$ ranges from 0 to 360 degrees in a full cycle.

The physical domain of the problem is defined by means of a single bounded subdomain. The scaling center of the subdomain is located at its centroid, and each edge is split into four equal segments. The semi-infinite nature of the half-space is included using eight unbounded subdomains, four of which are defined along the lower edge, and the other four along the right edge of the physical domain. Different excitation frequencies, degrees of polynomials, and discretization steps are used to analyze the results in detail. The verification process begins with comparing the numerical solutions for a combination of an overkill radial discretization and high-degree polynomials against the analytical data. To this end, the radial lines along the local coordinate $\xi$ are divided into 100 intervals, and the field variables along the circumferential coordinate $\eta$ are interpolated by means of 10 Gauss--Lobatto points. Next, the sensitivity of the solutions with respect to the number of radial intervals and circumferential points are analyzed. In all cases, the starting point in the bounded and unbounded subdomains are chosen to be $10^{-6}$ and 2, respectively. The damping ratio of the unbounded subdomains varies linearly from 1.0 to 0.05 by moving from the starting point to the open end. Three cases of excitation frequency, including 15, 25, and 35 Hz are chosen, for which the results are presented in the following.

\subsection{Frequency of 15 Hz}
\begin{figure}
    \centering
    \includegraphics[scale=1.0]{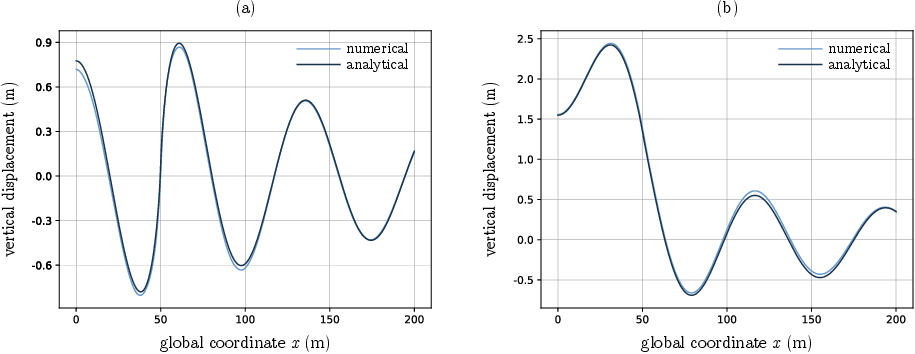}
    \caption{Vertical displacement of the surface for $f=15$ Hz: (a) real part, (b) imaginary part.}
    \label{fig:f15-curves}
\end{figure}

First, the solution method is assessed at a relatively low frequency of 15 Hz. To this end, the numerical result of the overkill model is plotted on the analytical data in Figure~\ref{fig:f15-curves}. The real and imaginary parts of the vertical displacement for the surface of the half-space are used for this comparison. As it can be seen, the curves obtained by the overkill numerical model almost completely lie on those of the experimental solution.

\begin{figure}
    \centering
    \includegraphics[scale=1.0]{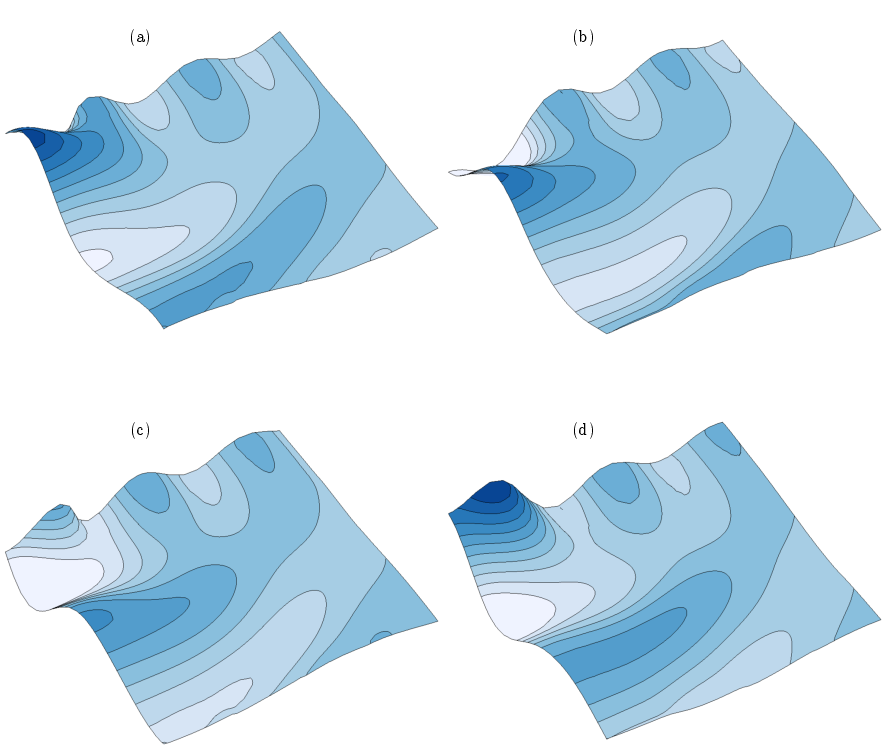}
    \caption{Vertical displacement field for $f=15$ Hz in different angles: (a) 0 degrees, (b) 90 degrees, (c) 180 degrees, (d) 270 degrees.}
    \label{fig:f15-contour}
\end{figure}

The vertical displacement field of the whole physical domain for different angles of 0, 90, 180, and 270 degrees are illustrated in Figure~\ref{fig:f15-contour}. Note that this solution field is reproduced within a single subdomain of the scaled boundary finite element method. The presented solution method enables the recovery of such highly-detailed field variables, which is absent in the available solution methods, such as the continued fraction--based methods in the works of Song and co-authors~\cite{song2007boundary, bazyar2008continued, song2009scaled, birk2012improved, chen2014high} and the shooting approach of Daneshyar et al.~\cite{daneshyar2021shooting}.

\begin{figure}
    \centering
    \includegraphics[scale=1.0]{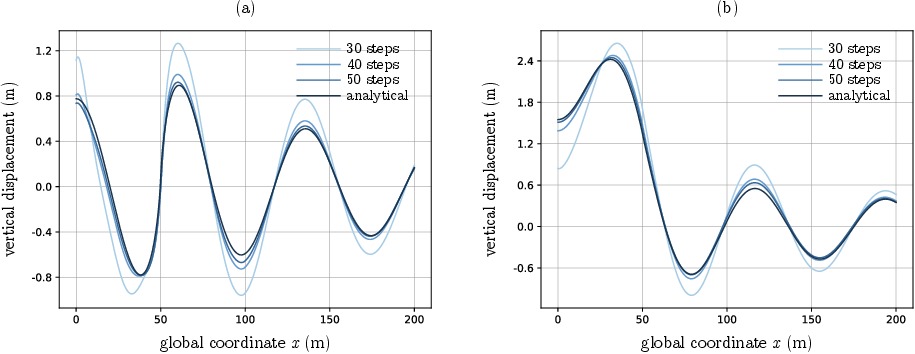}
    \caption{Vertical displacement of the surface for $f=15$ Hz using different radial discretization steps: (a) real part, (b) imaginary part.}
    \label{fig:f15-curves-steps}
\end{figure}

\begin{figure}
    \centering
    \includegraphics[scale=1.0]{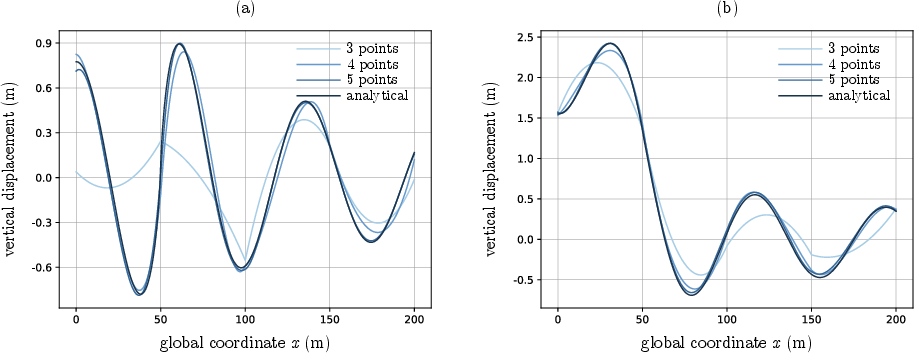}
    \caption{Vertical displacement of the surface for $f=15$ Hz using different number of interpolation points: (a) real part, (b) imaginary part.}
    \label{fig:f15-curves-points}
\end{figure}

\begin{figure}
    \centering
    \includegraphics[scale=1.0]{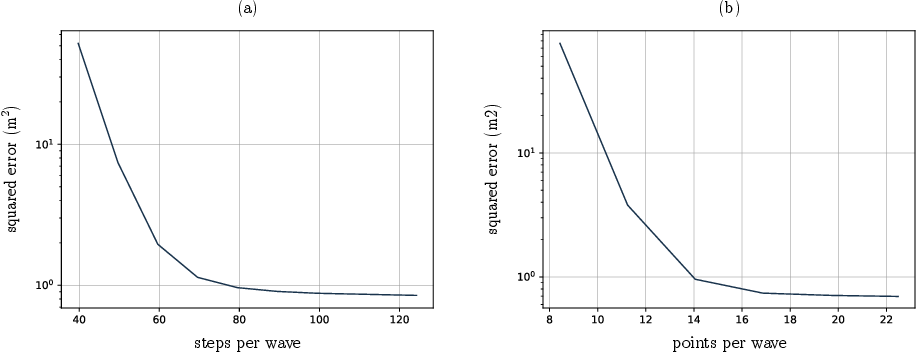}
    \caption{Sum of squared error for $f=15$ Hz: (a) versus the number of discretization steps per wave, (b) versus the number of interpolation points per wave.}
    \label{fig:f15-curves-errors}
\end{figure}

Now, the sensitivity of the results with respect to the discretization steps and interpolation functions is pursued. The first comparison proceeds by analyzing the real and imaginary parts of the vertical displacement of the surface for three cases of 30, 40, and 50 radial discretization steps in Figure~\ref{fig:f15-curves-steps}. As it was expected, the numerical solutions converge to the analytical one upon increasing the number of steps. Note that, to isolate the effect of radial discretization, all the cases are analyzed using 10 Gauss--Lobatto points. The same comparison is made for the degree of interpolation functions along the circumferential coordinate $\eta$. To this end, the number of radial discretization steps is kept fixed and only the number of interpolation points on the boundary elements of the subdomains are varied. Figure~\ref{fig:f15-curves-points} illustrates the accuracy of the real and imaginary parts of the results for the three cases of Gauss--Lobatto--Legendre polynomials interpolated by means of 3, 4, and 5 boundary points. A similar deduction to that of the cases with different radial discretization steps can be made. Lastly, the convergence trends upon increasing the number of discretization steps and degree of interpolation functions are depicted on logarithmic plots in Figure~\ref{fig:f15-curves-errors}. The vertical axes denotes the sum of the squared error of the vertical displacement of the surface involving both real and imaginary parts. Note that, to provide a rough estimation on the choice of discretization steps and interpolation functions, the horizontal axes are normalized with respect to the geometry of the subdomain and the frequency of the propagating wave so that they show how many discretization steps or interpolation points are defined along the wavelength
\begin{equation}
    \lambda_p = \frac{c_p}{f},
\end{equation}
wherein $c_p$ is the pressure wave speed and $f$ is the frequency of the propagating wave. Since all the radial distances of a subdomain are discretized by means of equal steps, the longest radial distance is critical and leads to the largest step. On the other hand, the longest boundary element of the subdomain is used to normalize the horizontal axis in the second type of the sensitivity analysis. Referring to the plots in Figure~\ref{fig:f15-curves-errors}, the numerical results reach a plateau beyond approximately 90 discretization steps per pressure wave for the first type of analysis and 17 interpolation points per pressure wave for the second type.
 
\subsection{Frequency of 25 Hz}
\begin{figure}
    \centering
    \includegraphics[scale=1.0]{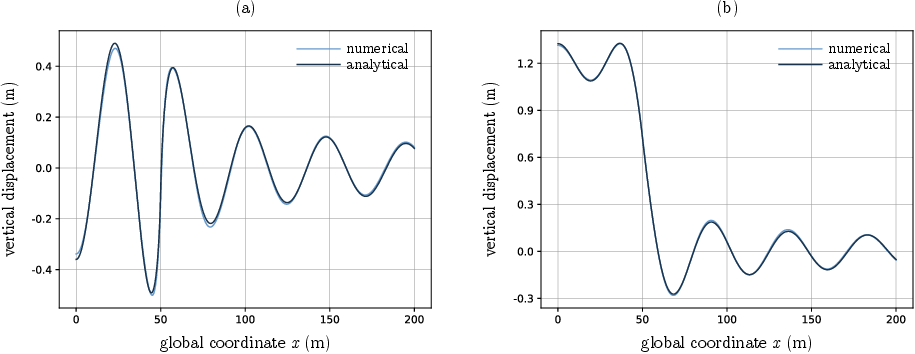}
    \caption{Vertical displacement of the surface for $f=25$ Hz: (a) real part, (b) imaginary part.}
    \label{fig:f25-curves}
\end{figure}

\begin{figure}
    \centering
    \includegraphics[scale=1.0]{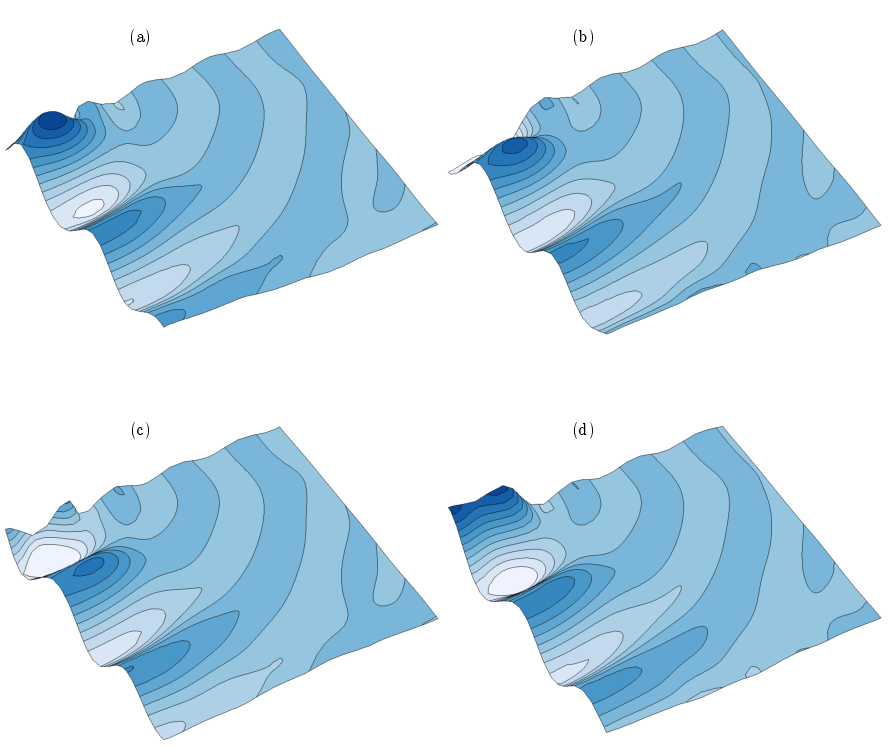}
    \caption{Vertical displacement field for $f=25$ Hz in different angles: (a) 0 degrees, (b) 90 degrees, (c) 180 degrees, (d) 270 degrees.}
    \label{fig:f25-contour}
\end{figure}

Here, the results of the presented method is analyzed for the moderate frequency 25 Hz. Similar to the previous comparison, the overkill numerical model defined by means of 10 Gauss--Lobatto points is approximated using 100 radial discretization steps. The robustness of the solution method can be simply deduced by comparing the real and imaginary parts of the vertical displacement of the surface of the half-space with the analytical curves in Figure~\ref{fig:f25-curves}. More importantly, it provides the inter-subdomain solution variables such as the vertical displacement field shown in Figure~\ref{fig:f25-contour}. It is worth to restate that this feature enables deploying the full potential of the scaled boundary finite element method.

\begin{figure}
    \centering
    \includegraphics[scale=1.0]{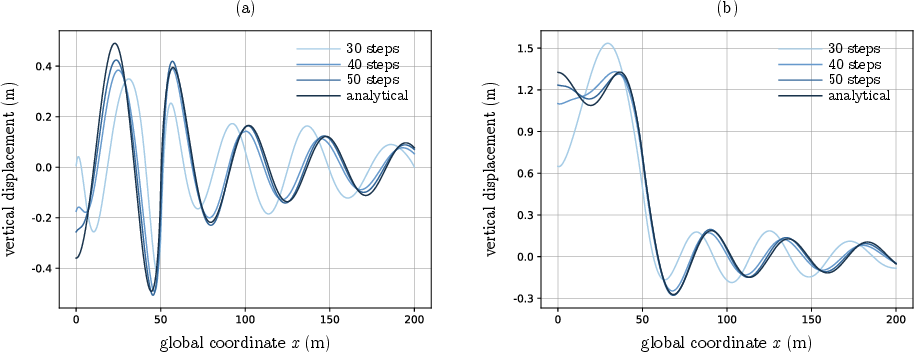}
    \caption{Vertical displacement of the surface for $f=25$ Hz using different radial discretization steps: (a) real part, (b) imaginary part.}
    \label{fig:f25-curves-steps}
\end{figure}

\begin{figure}
    \centering
    \includegraphics[scale=1.0]{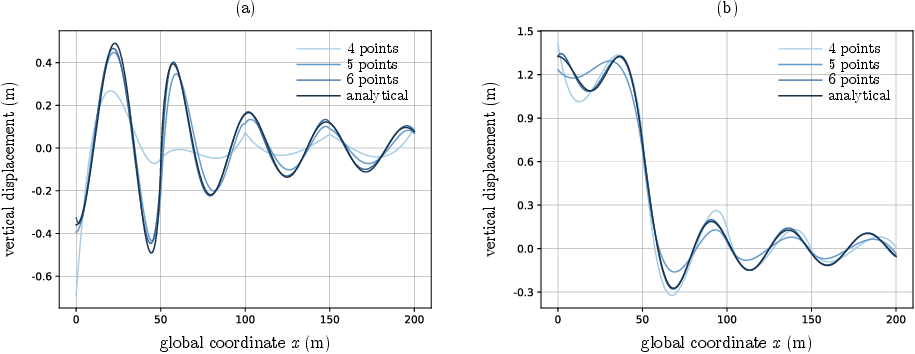}
    \caption{Vertical displacement of the surface for $f=25$ Hz using different number of interpolation points: (a) real part, (b) imaginary part.}
    \label{fig:f25-curves-points}
\end{figure}

\begin{figure}
    \centering
    \includegraphics[scale=1.0]{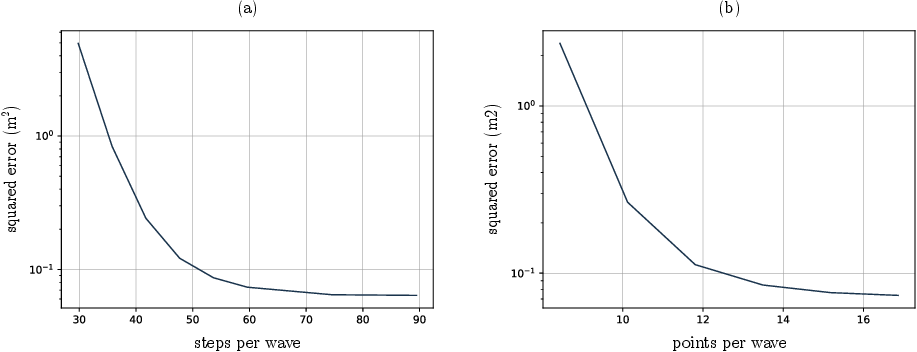}
    \caption{Sum of squared error for $f=25$ Hz: (a) versus the number of discretization steps per wave, (b) versus the number of interpolation points per wave.}
    \label{fig:f25-curves-errors}
\end{figure}

Next, analyzing the convergence trend of the numerical results is aimed at once again. For this purpose, the quality of approximations are assessed through varying the number of radial discretization steps and the degree of Gauss--Lobatto--Legendre polynomials. The former is carried out for 30, 40, and 50 discretization steps and the latter for 4, 5, and 6 Gauss--Lobatto points. Figure~\ref{fig:f25-curves-steps} shows the results of the cases with different radial discretization steps while Figure~\ref{fig:f25-curves-points} illustrates those of the cases with different number of interpolation points. The logarithmic error plots are also presented in Figure~\ref{fig:f25-curves-errors}. The first type of analysis enters into its plateau phase after approximately 75 steps, while the second one reaches that stage by passing 15 interpolation points.

\subsection{Frequency of 35 Hz}
\begin{figure}
    \centering
    \includegraphics[scale=1.0]{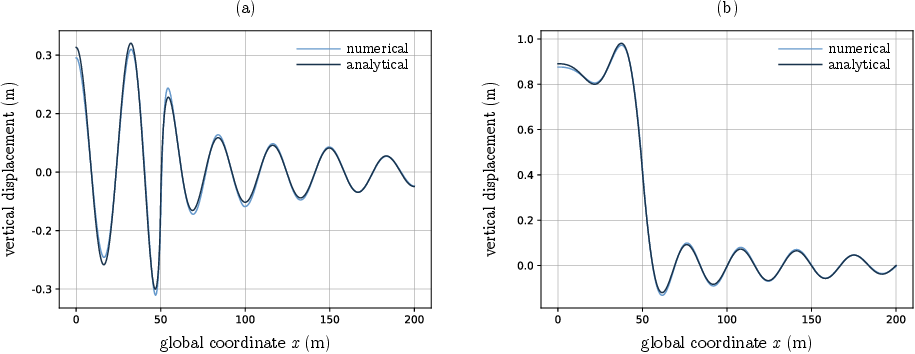}
    \caption{Vertical displacement of the surface for $f=35$ Hz: (a) real part, (b) imaginary part.}
    \label{fig:f35-curves}
\end{figure}

Lastly, the relatively high frequency case of 35 Hz is chosen to evaluate the solution method. The assessment proceeds similar to the previous cases by analyzing the overkill case. To this end, the numerical model of the half-space defined using 10 interpolation points along each boundary segments is approximated by means of 100 radial discretization steps. The real and imaginary parts of the vertical displacement of the surface for this model are plotted on their respective analytical curve in Figure~\ref{fig:f35-curves}. The close agreement between the numerical and analytical results demonstrates the robustness of the solution procedure once again.

\begin{figure}
    \centering
    \includegraphics[scale=1.0]{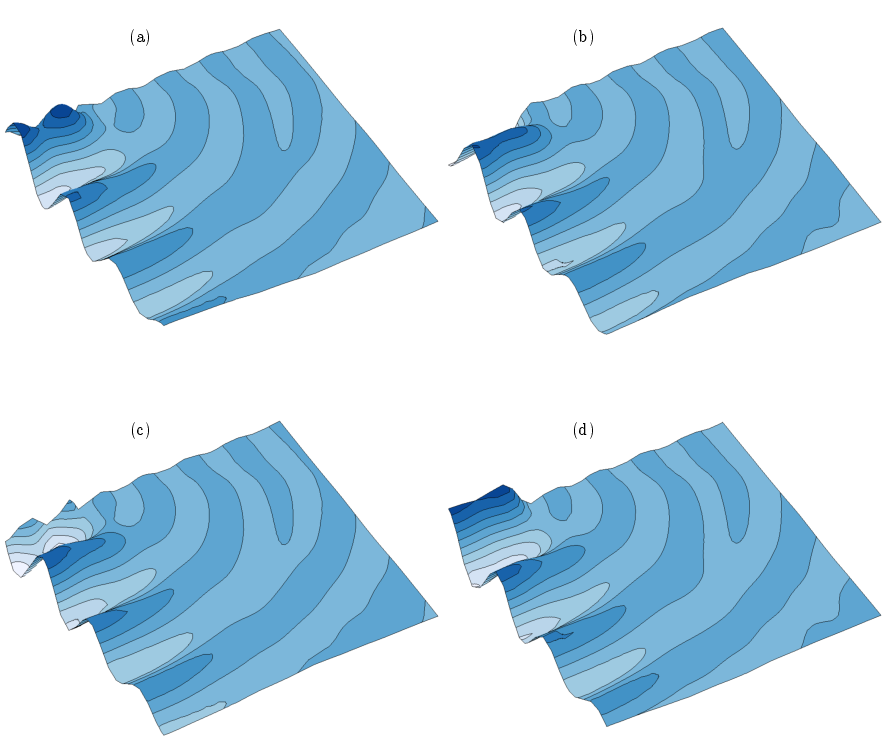}
    \caption{Vertical displacement field for $f=35$ Hz in different angles: (a) 0 degrees, (b) 90 degrees, (c) 180 degrees, (d) 270 degrees.}
    \label{fig:f35-contour}
\end{figure}

\begin{figure}
    \centering
    \includegraphics[scale=1.0]{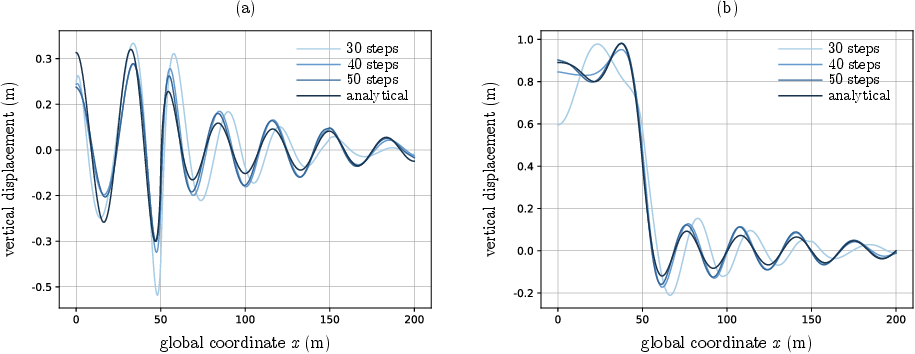}
    \caption{Vertical displacement of the surface for $f=35$ Hz using different radial discretization steps: (a) real part, (b) imaginary part.}
    \label{fig:f35-curves-steps}
\end{figure}

\begin{figure}
    \centering
    \includegraphics[scale=1.0]{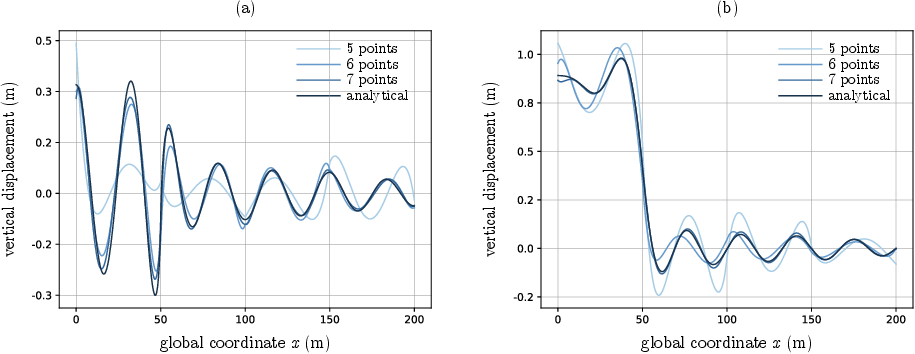}
    \caption{Vertical displacement of the surface for $f=35$ Hz using different number of interpolation points: (a) real part, (b) imaginary part.}
    \label{fig:f35-curves-points}
\end{figure}

\begin{figure}
    \centering
    \includegraphics[scale=1.0]{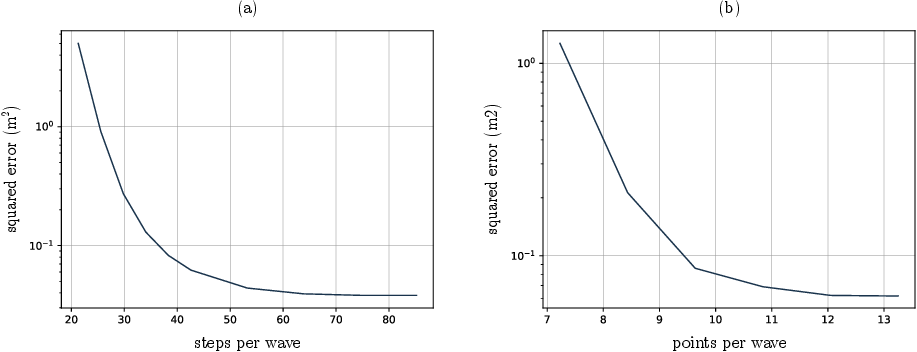}
    \caption{Sum of squared error for $f=35$ Hz: (a) versus the number of discretization steps per wave, (b) versus the number of interpolation points per wave.}
    \label{fig:f35-curves-errors}
\end{figure}

Next, we assess the convergence trends with respect to the number of steps and interpolation points. Three cases of 30, 40, and 50 steps and 5, 6, and 7 interpolation points per pressure wave are analyzed, for which the curves are compared against the analytical data in Figure~\ref{fig:f35-curves-steps} and \ref{fig:f35-curves-points}, respectively. The logarithmic error plots are also presented in Figure~\ref{fig:f35-curves-errors}. The plateau phase is reached after 65 steps per pressure wave and 12 interpolation points, respectively.

\section{Conclusion} \label{sec:conclusion}
We approached the scaled boundary finite element equations of elastodynamic problems in the frequency domain by means of a new solution technique. Upon discretizing the local radial directions, the dynamic stiffness matrix of bounded and unbounded subdomains can be directly extracted from the original second-order differential equations system. This technique enables the recovery of solution fields over the scaled boundary finite element spaces. As a result, the superiority of the method regarding its capability in reproducing highly-detailed solution fields over large subdomains can be put into practice. The task is achieved without introducing any new unknowns to the global algebraic system of equations. The analytical solution of a viscoelastic half-space subjected to harmonic strip loading was used to verify the solution method. The numerical results were in perfect agreement with the analytical ones. The highly-detailed inner-subdomain solution fields were also presented, showing the superiority of the presented solution method over its predecessors. Sensitivity analyses regarding the radial discretization and interpolation functions were also presented. The numerical solutions fall into their plateau phase after approximately 90, 75, and 65 steps per pressure wave for the frequencies 15, 25, and 35 Hz, respectively. The same occurred for 17, 15, and 12 interpolation points per pressure wave. Considering the fact that the wavelength of high-frequency waves are shorter than that of low-frequency waves, 65 discretization steps and 12 interpolation points per pressure wave would lead to asymptotic accuracy for all three cases. As a general rule of thumb, those choices can be made in accordance with the highest frequency involved.

\section{Declaration of competing interest} \label{sec:declaration}
The authors declare that they have no known competing financial interests or personal relationships that could have appeared to influence the work reported in this paper.

\section{Acknowledgement} \label{sec:acknowledgement}
This work was supported by the Alexander von Humboldt and Carl Friedrich von Siemens Foundations.

\bibliographystyle{ieeetr}
\bibliography{preprint}

\begin{thebibliography}{10}

\bibitem{song1997scaled}
C.~Song and J.~P. Wolf, ``The scaled boundary finite-element method—alias consistent infinitesimal finite-element cell method—for elastodynamics,'' {\em Computer Methods in applied mechanics and engineering}, vol.~147, no.~3-4, pp.~329--355, 1997.

\bibitem{reichel2023non}
R.~Reichel and S.~Klinkel, ``A non-uniform rational b-splines enhanced finite element formulation based on the scaled boundary parameterization for the analysis of heterogeneous solids,'' {\em International Journal for Numerical Methods in Engineering}, vol.~124, no.~9, pp.~2068--2092, 2023.

\bibitem{gravenkamp2018scaled}
H.~Gravenkamp and S.~Natarajan, ``Scaled boundary polygons for linear elastodynamics,'' {\em Computer Methods in Applied Mechanics and Engineering}, vol.~333, pp.~238--256, 2018.

\bibitem{gravenkamp2020mass}
H.~Gravenkamp, C.~Song, and J.~Zhang, ``On mass lumping and explicit dynamics in the scaled boundary finite element method,'' {\em Computer Methods in Applied Mechanics and Engineering}, vol.~370, p.~113274, 2020.

\bibitem{gravenkamp2020three}
H.~Gravenkamp, A.~A. Saputra, and S.~Eisentr{\"a}ger, ``Three-dimensional image-based modeling by combining sbfem and transfinite element shape functions,'' {\em Computational Mechanics}, vol.~66, no.~4, pp.~911--930, 2020.

\bibitem{zhang2021massively}
J.~Zhang, A.~Ankit, H.~Gravenkamp, S.~Eisentr{\"a}ger, and C.~Song, ``A massively parallel explicit solver for elasto-dynamic problems exploiting octree meshes,'' {\em Computer Methods in Applied Mechanics and Engineering}, vol.~380, p.~113811, 2021.

\bibitem{qu2022time}
Y.~Qu, J.~Zhang, S.~Eisentr{\"a}ger, and C.~Song, ``A time-domain approach for the simulation of three-dimensional seismic wave propagation using the scaled boundary finite element method,'' {\em Soil Dynamics and Earthquake Engineering}, vol.~152, p.~107011, 2022.

\bibitem{zhang2022asynchronous}
J.~Zhang, M.~Zhao, S.~Eisentr{\"a}ger, X.~Du, and C.~Song, ``An asynchronous parallel explicit solver based on scaled boundary finite element method using octree meshes,'' {\em Computer Methods in Applied Mechanics and Engineering}, vol.~401, p.~115653, 2022.

\bibitem{eisentrager2020sbfem}
J.~Eisentr{\"a}ger, J.~Zhang, C.~Song, and S.~Eisentr{\"a}ger, ``An sbfem approach for rate-dependent inelasticity with application to image-based analysis,'' {\em International Journal of Mechanical Sciences}, vol.~182, p.~105778, 2020.

\bibitem{chasapi2020geometrically}
M.~Chasapi and S.~Klinkel, ``Geometrically nonlinear analysis of solids using an isogeometric formulation in boundary representation,'' {\em Computational Mechanics}, vol.~65, pp.~355--373, 2020.

\bibitem{chasapi2020patch}
M.~Chasapi, W.~Dornisch, and S.~Klinkel, ``Patch coupling in isogeometric analysis of solids in boundary representation using a mortar approach,'' {\em International Journal for Numerical Methods in Engineering}, vol.~121, no.~14, pp.~3206--3226, 2020.

\bibitem{chasapi2021isogeometric}
M.~Chasapi, L.~Mester, B.~Simeon, and S.~Klinkel, ``Isogeometric analysis of 3d solids in boundary representation for problems in nonlinear solid mechanics and structural dynamics,'' {\em International Journal for Numerical Methods in Engineering}, 2021.

\bibitem{wallner2020scaled}
M.~Wallner, C.~Birk, and H.~Gravenkamp, ``A scaled boundary finite element approach for shell analysis,'' {\em Computer Methods in Applied Mechanics and Engineering}, vol.~361, p.~112807, 2020.

\bibitem{li2020efficient}
J.~Li, Z.~Shi, L.~Liu, and C.~Song, ``An efficient scaled boundary finite element method for transient vibro-acoustic analysis of plates and shells,'' {\em Computers \& Structures}, vol.~231, p.~106211, 2020.

\bibitem{klassen2021isogeometric}
M.~Klassen and S.~Klinkel, ``An isogeometric scaled boundary plate formulation for the analysis of ionic electroactive paper,'' {\em Acta Mechanica}, pp.~1--13, 2021.

\bibitem{song2002semi}
C.~Song and J.~P. Wolf, ``Semi-analytical representation of stress singularities as occurring in cracks in anisotropic multi-materials with the scaled boundary finite-element method,'' {\em Computers \& Structures}, vol.~80, no.~2, pp.~183--197, 2002.

\bibitem{yang2006fully}
Z.~Yang, ``Fully automatic modelling of mixed-mode crack propagation using scaled boundary finite element method,'' {\em Engineering Fracture Mechanics}, vol.~73, no.~12, pp.~1711--1731, 2006.

\bibitem{bulling2019high}
J.~Bulling, H.~Gravenkamp, and C.~Birk, ``A high-order finite element technique with automatic treatment of stress singularities by semi-analytical enrichment,'' {\em Computer Methods in Applied Mechanics and Engineering}, vol.~355, pp.~135--156, 2019.

\bibitem{pramod2019adaptive}
A.~Pramod, R.~Annabattula, E.~Ooi, C.~Song, S.~Natarajan, {\em et~al.}, ``Adaptive phase-field modeling of brittle fracture using the scaled boundary finite element method,'' {\em Computer Methods in Applied Mechanics and Engineering}, vol.~355, pp.~284--307, 2019.

\bibitem{zhang2020adaptive}
J.~Zhang, S.~Natarajan, E.~T. Ooi, and C.~Song, ``Adaptive analysis using scaled boundary finite element method in 3d,'' {\em Computer Methods in Applied Mechanics and Engineering}, vol.~372, p.~113374, 2020.

\bibitem{ooi2020polygon}
E.~Ooi, M.~Iqbal, C.~Birk, S.~Natarajan, E.~Ooi, and C.~Song, ``A polygon scaled boundary finite element formulation for transient coupled thermoelastic fracture problems,'' {\em Engineering Fracture Mechanics}, vol.~240, p.~107300, 2020.

\bibitem{ya2021open}
S.~Ya, S.~Eisentr{\"a}ger, C.~Song, and J.~Li, ``An open-source abaqus implementation of the scaled boundary finite element method to study interfacial problems using polyhedral meshes,'' {\em Computer Methods in Applied Mechanics and Engineering}, vol.~381, p.~113766, 2021.

\bibitem{iqbal2021development}
M.~Iqbal, C.~Birk, E.~Ooi, and H.~Gravenkamp, ``Development of the scaled boundary finite element method for crack propagation modeling of elastic solids subjected to coupled thermo-mechanical loads,'' {\em Computer Methods in Applied Mechanics and Engineering}, vol.~387, p.~114106, 2021.

\bibitem{gravenkamp2017efficient}
H.~Gravenkamp, A.~A. Saputra, C.~Song, and C.~Birk, ``Efficient wave propagation simulation on quadtree meshes using sbfem with reduced modal basis,'' {\em International Journal for Numerical Methods in Engineering}, vol.~110, no.~12, pp.~1119--1141, 2017.

\bibitem{liu2019automatic}
L.~Liu, J.~Zhang, C.~Song, C.~Birk, and W.~Gao, ``An automatic approach for the acoustic analysis of three-dimensional bounded and unbounded domains by scaled boundary finite element method,'' {\em International Journal of Mechanical Sciences}, vol.~151, pp.~563--581, 2019.

\bibitem{daneshyar2021general}
A.~Daneshyar and M.~Ghaemian, ``A general solution procedure for the scaled boundary finite element method via shooting technique,'' {\em Computer Methods in Applied Mechanics and Engineering}, vol.~384, p.~113996, 2021.

\bibitem{song1998scaled}
C.~Song and J.~P. Wolf, ``The scaled boundary finite-element method: analytical solution in frequency domain,'' {\em Computer Methods in Applied Mechanics and Engineering}, vol.~164, no.~1-2, pp.~249--264, 1998.

\bibitem{yann2004coupling}
J.~Yann, C.~Zhang, and F.~Jin, ``A coupling procedure of fe and sbfe for soil--structure interaction in the time domain,'' {\em International Journal for Numerical Methods in Engineering}, vol.~59, no.~11, pp.~1453--1471, 2004.

\bibitem{paronesso1998recursive}
A.~Paronesso and J.~P. Wolf, ``Recursive evaluation of interaction forces and property matrices from unit-impulse response functions of unbounded medium based on balancing approximation,'' {\em Earthquake engineering \& structural dynamics}, vol.~27, no.~6, pp.~609--618, 1998.

\bibitem{genes2012dynamic}
M.~C. Genes, ``Dynamic analysis of large-scale ssi systems for layered unbounded media via a parallelized coupled finite-element/boundary-element/scaled boundary finite-element model,'' {\em Engineering Analysis with Boundary Elements}, vol.~36, no.~5, pp.~845--857, 2012.

\bibitem{schauer2012parallel}
M.~Schauer, J.~E. Roman, E.~S. Quintana-Ort{\'\i}, and S.~Langer, ``Parallel computation of 3-d soil-structure interaction in time domain with a coupled fem/sbfem approach,'' {\em Journal of Scientific Computing}, vol.~52, pp.~446--467, 2012.

\bibitem{yang2007frobenius}
Z.~Yang, A.~Deeks, and H.~Hao, ``A frobenius solution to the scaled boundary finite element equations in frequency domain for bounded media,'' {\em International journal for numerical methods in engineering}, vol.~70, no.~12, pp.~1387--1408, 2007.

\bibitem{song2007boundary}
C.~Song and M.~H. Bazyar, ``A boundary condition in pad{\'e} series for frequency-domain solution of wave propagation in unbounded domains,'' {\em International Journal for Numerical Methods in Engineering}, vol.~69, no.~11, pp.~2330--2358, 2007.

\bibitem{bazyar2008continued}
M.~H. Bazyar and C.~Song, ``A continued-fraction-based high-order transmitting boundary for wave propagation in unbounded domains of arbitrary geometry,'' {\em International Journal for Numerical Methods in Engineering}, vol.~74, no.~2, pp.~209--237, 2008.

\bibitem{song2009scaled}
C.~Song, ``The scaled boundary finite element method in structural dynamics,'' {\em International Journal for Numerical Methods in Engineering}, vol.~77, no.~8, pp.~1139--1171, 2009.

\bibitem{birk2012improved}
C.~Birk, S.~Prempramote, and C.~Song, ``An improved continued-fraction-based high-order transmitting boundary for time-domain analyses in unbounded domains,'' {\em International Journal for Numerical Methods in Engineering}, vol.~89, no.~3, pp.~269--298, 2012.

\bibitem{chen2014high}
D.~Chen, C.~Birk, C.~Song, and C.~Du, ``A high-order approach for modelling transient wave propagation problems using the scaled boundary finite element method,'' {\em International Journal for Numerical Methods in Engineering}, vol.~97, no.~13, pp.~937--959, 2014.

\bibitem{daneshyar2021shooting}
A.~Daneshyar, P.~Sotoudeh, and M.~Ghaemian, ``A shooting approach to the scaled boundary finite element equations of elastodynamics in the frequency domain,'' {\em Computer Methods in Applied Mechanics and Engineering}, vol.~387, p.~114170, 2021.

\bibitem{daneshyar2023scaled}
A.~Daneshyar, P.~Sotoudeh, and M.~Ghaemian, ``The scaled boundary finite element method for dispersive wave propagation in higher-order continua,'' {\em International Journal for Numerical Methods in Engineering}, vol.~124, no.~4, pp.~880--927, 2023.

\bibitem{vu2006use}
T.~H. Vu and A.~J. Deeks, ``Use of higher-order shape functions in the scaled boundary finite element method,'' {\em International Journal for Numerical Methods in Engineering}, vol.~65, no.~10, pp.~1714--1733, 2006.

\bibitem{riley1999mathematical}
K.~F. Riley, M.~P. Hobson, and S.~J. Bence, {\em Mathematical methods for physics and engineering}.
\newblock American Association of Physics Teachers, 1999.

\bibitem{wolf2003scaled}
J.~P. Wolf, {\em The scaled boundary finite element method}.
\newblock John Wiley \& Sons, 2003.

\bibitem{song2004matrix}
C.~Song, ``A matrix function solution for the scaled boundary finite-element equation in statics,'' {\em Computer Methods in Applied Mechanics and Engineering}, vol.~193, no.~23-26, pp.~2325--2356, 2004.

\bibitem{song2018scaled}
C.~Song, {\em The Scaled Boundary Finite Element Method: Introduction to Theory and Implementation}.
\newblock Wiley, 2018.

\bibitem{ooi2016construction}
E.~T. Ooi, C.~Song, and S.~Natarajan, ``Construction of high-order complete scaled boundary shape functions over arbitrary polygons with bubble functions,'' {\em International Journal for Numerical Methods in Engineering}, vol.~108, no.~9, pp.~1086--1120, 2016.

\bibitem{song2006dynamic}
C.~Song, ``Dynamic analysis of unbounded domains by a reduced set of base functions,'' {\em Computer Methods in Applied Mechanics and Engineering}, vol.~195, no.~33-36, pp.~4075--4094, 2006.

\bibitem{kausel2019numerical}
E.~Kausel and H.~Gravenkamp, ``On the numerical solution of matrix bessel equations,'' {\em ZAMM-Journal of Applied Mathematics and Mechanics/Zeitschrift f{\"u}r Angewandte Mathematik und Mechanik}, vol.~99, no.~8, p.~e201800288, 2019.

\bibitem{patankar2018numerical}
S.~Patankar, {\em Numerical heat transfer and fluid flow}.
\newblock Taylor \& Francis, 2018.

\bibitem{miller1954field}
G.~Miller and H.~Pursey, ``The field and radiation impedance of mechanical radiators on the free surface of a semi-infinite isotropic solid,'' {\em Proceedings of the Royal Society of London. Series A. Mathematical and Physical Sciences}, vol.~223, no.~1155, pp.~521--541, 1954.

\end{thebibliography}

\end{document}